\documentclass[reqno]{amsart}

\usepackage{amsmath,amssymb,amscd,latexsym,amsthm}

\theoremstyle{theorem}
\newtheorem{theorem}{Theorem}

\theoremstyle{theorem}
\newtheorem{corollary}{Corollary}

\theoremstyle{definition}

\begin{document}

\title[Solutions of equations over finite fields]{Solutions of equations over finite fields: enumeration via bijections}

\author{Ioulia N. Baoulina}

\address{Department of Mathematics, Moscow State Pedagogical University, Krasnoprudnaya str. 14, Moscow 107140, Russia}
\email{jbaulina@mail.ru}

\date{}

\maketitle

\begin{abstract}
We present a simple proof of the well-known fact concerning the number of solutions of diagonal equations over finite fields. In a similar manner, we give an alternative proof of the recent result on generalizations of Carlitz equations. In both cases, the use of character sums is avoided by using an elementary combinatorial argument.
\end{abstract}

\section{Introduction.}

Let $\mathbb F_q$ be a finite field of $q$ elements and let $\mathbb F_q^*=\mathbb F_q\setminus\{0\}$. We consider equations of the form $f_1(x_1,\dots,x_n)=f_2(x_1,\dots,x_n)$, where $f_1$ and $f_2$ are polynomials with coefficients in $\mathbb F_q$. Explicit formulas for the number of solutions $(x_1,\dots,x_n)\in\mathbb F_q^n$ are known only for certain classes of equations (see \cite{BEW,LN,MP}). The traditional approach is to express the number of solutions in terms of character sums and to apply known results about these sums. In this note, we present an alternative approach, based on the observation that for some families of equations the number of solutions does not depend on the coefficients, that is, the same formula is valid for the whole family. This suggests that one should try to establish a bijective correspondence between the sets of solutions of equations belonging to a given family and to use this correspondence to obtain an explicit formula for the number of solutions. In Sections~\ref{s2} and \ref{s3}, we illustrate these ideas. We denote by $N[f_1(x_1,\dots,x_n)=f_2(x_1,\dots,x_n)]$ the number of solutions of the equation $f_1(x_1,\dots,x_n)=f_2(x_1,\dots,x_n)$ in $\mathbb F_q^n$ and by $N^*[f_1(x_1,\dots,x_n)=f_2(x_1,\dots,x_n)]$ the number of such solutions with ${x_1\cdots x_n\ne 0}$. Unless explicitly stated otherwise, we will assume that equations have coefficients in $\mathbb F_q^*$, exponents occurring in them are positive integers and $n\ge 2$. We write $|\mathcal{S}|$ for the number of elements of a finite set~$\mathcal{S}$.

\section{Diagonal equations.}\label{s2} 
We consider an equation of the type
$$
a_1^{}x_1^{m_1}+\dots+a_n^{}x_n^{m_n}=0.
$$
Such an equation is called a diagonal equation. As $x_j$ runs through all elements of $\mathbb F_q$, $x_j^{m_j}$ runs through the same elements as $x_j^{\gcd(m_j,q-1)}$ does with the same multiplicity. Therefore,
$$
N[a_1^{}x_1^{m_1}+\dots+a_n^{}x_n^{m_n}=0]=N[a_1^{}x_1^{d_1}+\dots+a_n^{}x_n^{d_n}=0],
$$
where $d_j=\gcd(m_j,q-1)$. In 1971, Joly~\cite{joly} obtained the following result.
\begin{theorem}
\label{t1}
Assume that $\gcd(d_j,d_1\cdots d_n/d_j)=1$ for some $j$. Then
$$
N[a_1^{}x_1^{d_1}+\dots+a_n^{}x_n^{d_n}=0]=q^{n-1}.
$$
\end{theorem}
In his proof, Joly used the expression for the number of solutions in terms of Gauss sums. We give an alternative proof based on the ideas mentioned in the introduction.
\begin{proof}
Without loss of generality we may assume that $\gcd(d_1,d_2\cdots d_n)=1$. For $c\in\mathbb F_q$, let $\mathcal{S}_c$ denote the set of solutions $(x_1,\dots,x_n)\in\mathbb F_q^n$ of the equation
$$
a_1^{}cx_1^{d_1}+a_2^{}x_2^{d_2}+\dots+a_n^{}x_n^{d_n}=0
$$
and $\mathcal{S}'_c$ denote the set of such solutions with $x_1\ne 0$. Clearly, 
$\mathcal{S}_{c_1}\setminus\mathcal{S}'_{c_1}=\mathcal{S}_{c_2}\setminus\mathcal{S}'_{c_2}$ for any $c_1,c_2\in\mathbb F_q$. In particular, we have $|\mathcal{S}_1\setminus\mathcal{S}'_1|=|\mathcal{S}_0\setminus\mathcal{S}'_0|$, or, equivalently,
\begin{equation}
\label{eq1}
|\mathcal{S}_1|=|\mathcal{S}'_1|+|\mathcal{S}_0|-|\mathcal{S}'_0|.
\end{equation}
Observe that for any $n$-tuple $(x_1,\dots,x_n)\in\mathbb F_q$ with $x_1\ne 0$ there exists a unique $c\in\mathbb F_q$ such that $(x_1,\dots,x_n)\in\mathcal{S}'_c$. In other words, we have
$$
\bigvee_{c\in\mathbb F_q}\mathcal{S}'_c=\mathbb F_q^*\times\underbrace{\mathbb F_q\times\dots\times\mathbb F_q}_{n-1}.
$$
This yields
\begin{equation}
\label{eq2}
\sum_{c\in\mathbb F_q}|\mathcal{S}'_c|=q^{n-1}(q-1).
\end{equation}
Since $\gcd(d_1,d_2\cdots d_n)=1$, we can choose a positive integer $t$ with $t\equiv 0\!\pmod{d_1}$ and $t\equiv 1\pmod{d_2\cdots d_n}$. It is easy to see that for a fixed $c\in\mathbb F_q^*$ the map
$$
(x_1,x_2,\dots,x_n) \longmapsto (c^{t/d_1}x_1,c^{(t-1)/d_2}x_2,\dots,c^{(t-1)/d_n}x_n)
$$
is a bijection between $\mathcal{S}'_c$ and $\mathcal{S}'_1$. This implies that $|\mathcal{S}'_c|=|\mathcal{S}'_1|$ for any $c\in\mathbb F_q^*$, so that \eqref{eq2} can be rewritten as $|\mathcal{S}'_0|+(q-1)|\mathcal{S}'_1|=q^{n-1}(q-1)$, or, equivalently,
\begin{equation}
\label{eq3}
|\mathcal{S}'_1|=q^{n-1}-\frac 1{q-1}\,|\mathcal{S}'_0|.
\end{equation}
Finally, note that
\begin{align*}
|\mathcal{S}_0|&=qN[a_2^{}x_2^{d_2}+\dots+a_n^{}x_n^{d_n}=0],\\
|\mathcal{S}'_0|&=(q-1)N[a_2^{}x_2^{d_2}+\dots+a_n^{}x_n^{d_n}=0],
\end{align*}
and thus
\begin{equation}
\label{eq4}
|\mathcal{S}_0|=\frac q{q-1}\,|\mathcal{S}'_0|.
\end{equation}
Combining \eqref{eq1}, \eqref{eq3} and \eqref{eq4}, we find that
$$
N[a_1^{}x_1^{d_1}+\dots+a_n^{}x_n^{d_n}=0]=|\mathcal{S}_1|=q^{n-1},
$$
as desired.
\end{proof}

The inclusion-exclusion principle of combinatorics yields the following corollary, which will be useful in the sequel.
\begin{corollary}
\label{c1}
If $d_1,\dots,d_n$ are pairwise coprime then
$$
N^*[a_1^{}x_1^{d_1}+\dots+a_n^{}x_n^{d_n}=0]=\frac{(q-1)^n+(-1)^n(q-1)}q.
$$
\end{corollary}

\section{Generalizations of Carlitz equations.}\label{s3}
Carlitz~\cite{carlitz1} studied equations of the form
$$
(x_1+\dots+x_n)^2=bx_1\cdots x_n
$$
over finite fields of odd characteristic. In particular, he proved that
$$
N[(x_1+x_2+x_3)^2=bx_1x_2x_3]=q^2+1.
$$
In \cite{baoulina}, the above result of Carlitz was generalized as follows:
$$
N[(x_1^{m_1}+\dots+x_n^{m_n})^k=bx_1^{}\cdots x_n^{}]=q^{n-1}+(-1)^{n-1},
$$
provided that $\gcd(\sum_{j=1}^n (M/m_j)-kM,(q-1)/D)=1$ and $d_1,\dots,d_n$ are pairwise coprime, where $d_j=\gcd(m_j,q-1)$, $j=1,\dots,n$, $M=\text{lcm}[m_1,\dots,m_n]$ and $D=\text{lcm}[d_1,\dots,d_n]$ (here $q$ is not necessarily odd). See the references in~\cite{baoulina} for earlier results in this direction. 

A further generalization of the Carlitz equation was recently considered in~\cite{PZC}. By using the augmented degree matrix and Gauss sums, the authors established the following: if $\gcd(\sum_{j=1}^n (k_jm_1\cdots m_n/m_j)-km_1\cdots m_n,q-1)=1$, then
$$
N[(a_1^{}x_1^{m_1}+\dots+a_n^{}x_n^{m_n})^k=bx_1^{k_1}\cdots x_n^{k_n}]=q^{n-1}+(-1)^{n-1}.
$$
Let us remark that, with the notation as above, we have
$$
\sum_{j=1}^n \frac{k_jm_1\cdots m_n}{m_j}-km_1\cdots m_n=\frac{m_1\cdots m_n}M\,\biggl(\sum_{j=1}^n \frac{k_jM}{m_j}-kM\biggr).
$$
Therefore the condition $\gcd(\sum_{j=1}^n (k_jm_1\cdots m_n/m_j)-km_1\cdots m_n,q-1)=1$ is equivalent to the following two conditions:
$$
\gcd\biggl(\sum_{j=1}^n \frac{k_jM}{m_j}-kM,q-1\biggr)=1,\qquad 
\gcd\biggl(\frac{m_1\cdots m_n}M,q-1\biggr)=1.
$$
It is readily seen that $M\cdot\gcd(m_i,m_j)\mid m_1\cdots m_n$. Since $\gcd(d_i,d_j)\mid(q-1)$ and $\gcd(d_i,d_j)\mid\gcd(m_i,m_j)$, we conclude that $\gcd(m_1\cdots m_n/M,q-1)=1$ if and only if $d_1,\dots,d_n$ are pairwise coprime. Summarizing, we reformulate the main result of~\cite{PZC} in the following way: 
\begin{theorem}
\label{t2}
Assume that $\gcd(\sum_{j=1}^n (k_jM/m_j)-kM,q-1)=1$ and $d_1,\dots,d_n$ are pairwise coprime. Then
$$
N[(a_1^{}x_1^{m_1}+\dots+a_n^{}x_n^{m_n})^k=bx_1^{k_1}\cdots x_n^{k_n}]=q^{n-1}+(-1)^{n-1}.
$$
\end{theorem}
Our proof of Theorem~\ref{t2} uses an argument similar to the one used in the proof of Theorem~\ref{t1}.
\begin{proof}
For $c\in\mathbb F_q$, let $\mathcal{S}_c$ denote the set of solutions $(x_1,\dots,x_n)\in\mathbb F_q^n$ of the equation
$$
(a_1^{}x_1^{m_1}+\dots+a_n^{}x_n^{m_n})^k=bcx_1^{k_1}\cdots x_n^{k_n}
$$
and $\mathcal{S}^*_c$ denote the set of such solutions with $x_1\cdots x_n\ne 0$. Proceeding by the same type of argument as in the proof of Theorem~\ref{t1}, we find that
\begin{equation}
\label{eq5}
|\mathcal{S}_1|=|\mathcal{S}^*_1|+|\mathcal{S}_0|-|\mathcal{S}^*_0|
\end{equation}
and
\begin{equation}
\label{eq6}
\sum_{c\in\mathbb F_q}|\mathcal{S}^*_c|=(q-1)^n.
\end{equation}
Let $g$ be a generator of the cyclic group $\mathbb F_q^*$. Then $g^{\sum_{j=1}^n (k_jM/m_j)-kM}$ is also a generator of $\mathbb F_q^*$, in view of the fact that $\gcd(\sum_{j=1}^n (k_jM/m_j)-kM,q-1)=1$. For a fixed $c\in\mathbb F_q^*$, there is an integer $t$ such that $c=g^{t(\sum_{j=1}^n (k_jM/m_j)-kM)}$. It is readily seen that the map
$$
(x_1,\dots,x_n) \longmapsto (g^{tM/m_1}x_1,\dots,g^{tM/m_n}x_n)
$$
is a bijection between $\mathcal{S}^*_1$ and $\mathcal{S}^*_c$, and, exactly as in the proof of Theorem~\ref{t1}, we deduce from \eqref{eq6} that
$$
|\mathcal{S}^*_1|=(q-1)^{n-1}-\frac 1{q-1}\,|\mathcal{S}^*_0|.
$$
Substituting this into \eqref{eq5} and recalling that
\begin{align*}
|\mathcal{S}_0|&=N[a_1^{}x_1^{d_1}+\dots+a_n^{}x_n^{d_n}=0]=q^{n-1},\\
|\mathcal{S}^*_0|&=N^*[a_1^{}x_1^{d_1}+\dots+a_n^{}x_n^{d_n}=0]=\frac{(q-1)^n+(-1)^n(q-1)}q,
\end{align*}
by Theorem~\ref{t1} and Corollary~\ref{c1}, we obtain the result.
\end{proof}

Note that we used the condition that $d_1,\dots,d_n$ are pairwise coprime only for determining the number of solutions of the corresponding diagonal equation. By relaxing this condition and proceeding exactly as above, we establish the following result.
\begin{theorem}
\label{t3}
Assume that $\gcd(\sum_{j=1}^n (k_jM/m_j)-kM,q-1)=1$. Then
\begin{multline*}
N[(a_1^{}x_1^{m_1}+\dots+a_n^{}x_n^{m_n})^k=bx_1^{k_1}\cdots x_n^{k_n}]=(q-1)^{n-1}\\
+N[a_1^{}x_1^{d_1}+\dots+a_n^{}x_n^{d_n}=0]-\frac q{q-1}N^*[a_1^{}x_1^{d_1}+\dots+a_n^{}x_n^{d_n}=0].
\end{multline*}
\end{theorem}
Theorem~\ref{t3} allows us to determine $N[(a_1^{}x_1^{m_1}+\dots+a_n^{}x_n^{m_n})^k=bx_1^{k_1}\cdots x_n^{k_n}]$ explicitly in certain cases when the number of solutions of the corresponding diagonal equation is known. In particular, using \cite[Theorem~2]{SY} and the inclusion-exclusion principle, we can generalize Theorem~\ref{t2} as follows (for a similar result in the case $a_1=\dots=a_n=1$, $k_1=\dots=k_n=1$ we refer to \cite{baoulina}).
\begin{theorem}
\label{t4}
Assume that $\gcd(\sum_{j=1}^n (k_jM/m_j)-kM,q-1)=1$, $d_1,\dots,d_t$ are odd, $d_{t+1},\dots,d_n$ are even, $d_1,\dots,d_t,d_{t+1}/2,\dots,d_n/2$ are pairwise coprime, ${0\le t\le n}$. Then
\begin{align*}
N[(a_1^{}x_1^{m_1}&+\dots+a_n^{}x_n^{m_n})^k=bx_1^{k_1}\cdots x_n^{k_n}]=q^{n-1}+(-1)^{n-1}\\
&+(-1)^{n-1}\sum_{j=1}^{[(n-t)/2]}\eta((-1)^j)\sigma_{2j}(\eta(a_{t+1}),\dots,\eta(a_n))q^j\\
&+\begin{cases}
\eta((-1)^{n/2}a_1\cdots a_n)q^{(n-2)/2}(q-1)&\text{if $t=0$ and $n$ is even,}\\
0&\text{otherwise,}
\end{cases}
\end{align*}
where $\sigma_{2j}(z_1,\dots,z_{n-t})$ are the elementary symmetric polynomials and $\eta$ is the quadratic character on $\mathbb F_q$ ($\eta(x)=+1,-1,0$ according as $x$ is a square, a non-square or zero in $\mathbb F_q$).
\end{theorem}

\section{Concluding remarks.}
The method used in the proof of Theorem~\ref{t2} can be applied to a more general class of equations, namely,
$$
f(x_1,\dots,x_n)=bx_1^{k_1}\cdots x_n^{k_n},
$$
where $f\in\mathbb F_q[x_1,\dots,x_n]$ is a polynomial for which there exist positive integers $r,r_1,\dots,r_n$ such that $f(c^{r_1}x_1,\dots,c^{r_n}x_n)=c^rf(x_1,\dots,x_n)$ for every $c\in\mathbb F_q$ (such a polynomial is sometimes called a quasi-homogeneous polynomial). If, in addition, $\gcd(\sum_{j=1}^n k_jr_j-r,q-1)=1$, then we have the expression
\begin{multline*}
N[f(x_1,\dots,x_n)=bx_1^{k_1}\cdots x_n^{k_n}]=(q-1)^{n-1}\\
+N[f(x_1,\dots,x_n)=0]-\frac q{q-1}N^*[f(x_1,\dots,x_n)=0],
\end{multline*}
and so the number of solutions of $f(x_1,\dots,x_n)=bx_1^{k_1}\cdots x_n^{k_n}$ is easily evaluated once $N[f(x_1,\dots,x_n)=0]$ and $N^*[f(x_1,\dots,x_n)=0]$ are known.


\end{document}